\documentclass[leqno,hyperref]{book}

\usepackage{hyperref}

\newfont{\neweu}{eurm10 scaled 1000} 

\def\eu#1{\mbox{\neweu #1}}

\def\bibname{\centerline{\textbf{REFERENCES}}}
\def\thebibliography#1{\paragraph*{\uppercase{\bibname}}\list
{{\bf [\arabic{enumi}]}}{\settowidth\labelwidth{[#1]}\leftmargin\labelwidth
\advance\leftmargin\labelsep\usecounter{enumi}}
\def\newblock{\hskip .11em plus .33em minus .07em}
\sloppy\clubpenalty4000\widowpenalty4000
\sfcode`\.=1000\relax}

\def\Dj{D{\hspace{-.75em}\raisebox{.3ex}{-}\hspace{.4em}}}

\begin{document}

\thispagestyle{plain}

\section{}

\centerline{\large \boldmath \bf A NOTE ABOUT THE $\{K_{i}(z)\}_{i=1}^{\infty}$ FUNCTIONS}
\footnotetext{Research partially supported by the MNTRS, Serbia, Grant No. 144020.}

\vspace*{5.00 mm}

\centerline{\large \it Branko J. Male\v sevi\' c}

\vspace*{5.00 mm}

\begin{center}
\parbox{25.14 cc}{\scriptsize
\textbf{\boldmath
In the article \mbox{\bf \cite{Petojevic_06}}, \textsc{A. Petojevi\' c} verified useful
properties of the $K_{i}(z)$ functions which generalize Kurepa's \mbox{\bf \cite{Kurepa_71}}
left factorial function. In this note, we present simplified proofs of two of these results
and we answer the open question stated in \mbox{\bf \cite{Petojevic_06}}. Finally, we
discuss the differential transcendency of the $K_{i}(z)$ functions.}}
\end{center}

\bigskip

{\sc A. Petojevi\' c} \mbox{\bf \cite[\mbox{\rm p.$\,$3.}]{Petojevic_02}} considered the family of functions:
\begin{equation}
\label{Eq_1}
\label{Fun_vmsaz}
{}_{v}M_{m}(s;a,z)
=
\displaystyle\sum\limits_{k=1}^{v}{(-1)^{k-1} {z+m+1-k \choose m+1}
\mbox{$\cal L$}[s;{}_{2}F_{1}(a,k-z,m+2;1-t)]},
\end{equation}
for $\Re(z)\,\!>\!\,v\!-\!m\!-\!2$, where $v \!\in\! \mbox{\bf N}$ is a positive integer;
$m \!\in\! \{-1, 0, 1, 2, \ldots\}$ is an integer; $s,a,z$ are complex variables;
$\mbox{$\cal L$}[s;F(t)]$ is {\sc Laplace} transform and~${}_{2}F_{1}(a,b,c;x)$ is the hypergeometric function
$(|x| \!<\! 1)$. {\sc \Dj .$\,$Ku\-repa} has considered in the articles \mbox{\bf \cite[\mbox{\rm p.$\,$151.}]{Kurepa_71}}
and \mbox{\bf \cite[\mbox{\rm p.$\,$297.}]{Kurepa_73}} a complex function defined by the integral:
\begin{equation}
\label{Eq_2}
K(z)
=
\displaystyle\int\limits_{0}^{\infty}{e^{-t} \displaystyle\frac{t^z-1}{t-1} \, dt},
\end{equation}
for \mbox{$\Re(z) \!>\! 0$}. Especially, for {\sc Ku\-repa}'s function  $K(z)$, it is true that
$K(z) \!=\! {}_{1}M_{0}(1;1,z)$, for~\mbox{$\Re(z) \!>\! 0$}, according to \mbox{\bf \cite{Petojevic_06}}.
For various of values of parameters $v,m,s,a,z$ from~(\ref{Fun_vmsaz}), different special functions,
as presented in \mbox{\bf \cite{Petojevic_06}}, are obtained. {\sc A. Petojevi\' c} has considered
in the article \mbox{\bf \cite[\mbox{\rm p.$\,$1640.}]{Petojevic_06}} the following sequence of functions:
\begin{equation}
\label{Eq_3}
K_{i}(z)
=
\displaystyle
\frac{{}_{1}M_{0}(1;1,z+i-1)-{}_{1}M_{0}(1;1,i-1)}{{}_{1}M_{-1}(1;1,i)},
\end{equation}
for $i \!\in\! \mbox{\bf N}$ and $\Re(z) \!>\! -i$. On the basis of the definition in (\ref{Eq_3}),
the following representation via {\sc Ku\-repa}'s function is true:
\begin{equation}
\label{Eq_4}
K_{i}(z)
=
\displaystyle\frac{1}{(i\!-\!1)!}
{\Big (}K(z+i-1) - K(i-1){\Big )},
\end{equation}
for $i \!\in\! \mbox{\bf N}$ and $\Re(z) \!>\! -i+1$. Note that $K(0) \!=\! 0$
\mbox{\bf \cite[\mbox{\rm p.$\,$297.}]{Kurepa_73}} and therefore $K_{1}(z) \!=\! K(z)$
for \mbox{$\Re(z) \!>\! 0$}. Analytical and differential--algebraic pro\-perties of {\sc Ku\-repa}'s
function $K(z)$ are considered in articles \mbox{\boldmath $[1 - 12]$} and in many other articles.
{\sl On the basis of well-known statements for {\sc Ku\-repa}'s function $K(z)$, using representation
{\rm $(\ref{Eq_4})$}, in many cases we can get simple proofs for analogous statements for $K_{i}(z)$ functions}.
For example, it is a well-known fact that it is possible to analytically continue {\sc Ku\-repa}'s function
to a meromorphic function with simple poles at integer points $z=-1$ and $z = -m$, ($m \geq 3$)
\mbox{\bf \cite[\mbox{\rm p.$\,$303.}]{Kurepa_73}}, \mbox{\bf \cite[\mbox{\rm p.$\,$474.}]{Slavic_73}}.
Residues of {\sc Ku\-repa}'s function at these poles have the~following~form~\mbox{\bf \cite{Kurepa_73}}:

\break

\noindent
\begin{equation}
\label{Eq_5}
\mathop{\mbox{\rm res}}\limits_{\mbox{\tiny $z = -1$}}{K(z)} = -1
\;\quad\;\mbox{\rm and}\quad
\mathop{\mbox{\rm res}}\limits_{\mbox{\tiny $z = -m$}}{K(z)} =
\displaystyle\sum\limits_{k=2}^{m-1}{\displaystyle\frac{(-1)^{k-1}}{k!}},
\;\; (m\!\geq\!3).
\end{equation}
For {\sc Ku\-repa}'s function $K(z)$ the infinite point is an essential singularity \mbox{\bf \cite{Slavic_73}}.
Hence, on the basis of (\ref{Eq_4}), each function $K_{i}(z)$ is meromorphic with simple poles at  integer points
$z\!=\!-i\,$ and $\,z\!=\!-(i+m)$, $(m \!\geq\! 2)$. On~the~basis of (\ref{Eq_4}) we~have:
\begin{equation}
\label{Eq_6}
\;\mathop{\mbox{\rm res}}\limits_{\mbox{\tiny $z=-(i\!+\!m)$}}{\!\!\!\!K_{i}(z)}
=
\displaystyle\frac{1}{(i\!-\!1)!}
\;\; \cdot \!\!\!\!
\mathop{\mbox{\rm res}}\limits_{\mbox{\tiny $z=-(i\!+\!m)$}}{\!\!\!\!K(z+i-1)}
=
\displaystyle\frac{1}{(i\!-\!1)!}
\;\; \cdot \!\!\!\!
\mathop{\mbox{\rm res}}\limits_{\mbox{\tiny $z=-(m\!+\!1)$}}{\!\!K(z)},
\end{equation}
where $m=0$ or $m \!\geq\! 2$. Hence:
\begin{equation}
\label{Eq_7}
\;\;\;\mathop{\mbox{\rm res}}\limits_{\mbox{\tiny $z=-i$}}{\!\!K_{i}(z)}
\!=\!
-\displaystyle\frac{1}{(i\!-\!1)!}
\;\;\quad\mbox{\rm and}\quad \!\!\!\!
\mathop{\mbox{\rm res}}\limits_{\mbox{\tiny $z=-(i\!+\!m)$}}{\!\!K_{i}(z)}
\!=\!
\displaystyle\frac{1}{(i\!-\!1)!}
\cdot\!\! \displaystyle\sum_{k=2}^{m}{
\displaystyle\frac{(-1)^{k-1}}{k!}}, \;\; (m\!\geq\!2).
\end{equation}
For each $K_{i}(z)$ function the infinite point is an essential singularity. Therefore, we get Theorem {\bf 3.3.}
from \mbox{\bf \cite{Petojevic_06}}. Next, it is a well-known fact that for {\sc Ku\-repa}'s function the following
asymptotic relation $K(x) \sim \Gamma(x)$ is true for real $x$ such that $x \rightarrow \infty$ and where $\Gamma(x)$
is the gamma function \mbox{\bf \cite[\mbox{\rm p.$\,$299.}]{Kurepa_73}}. Hence, for fixed $i \!\in\! \mbox{\bf N}$ and
real $x \!>\! -i\!+\!1$, on the basis of (\ref{Eq_4}),~we~get:
\begin{equation}
\label{Eq_8}
\displaystyle\frac{K_{i}(x)}{\Gamma(x+i-1)}
=
\displaystyle\frac{1}{(i\!-\!1)!}
\cdot
\displaystyle\frac{K(i+x-1)-K(i-1)}{\Gamma(x+i-1)}
\mathop{\longrightarrow}\limits_{\mbox{\tiny $x \!\rightarrow\! \infty$}}
\displaystyle\frac{1}{(i\!-\!1)!}
\end{equation}
and
\begin{equation}
\label{Eq_9}
\displaystyle\frac{K_{i}(x)}{\Gamma(x+i)}
=
\displaystyle\frac{1}{(i\!-\!1)!}
\cdot
\displaystyle\frac{K(i+x-1)-K(i-1)}{(x+i-1) \Gamma(x+i-1)}
\mathop{\longrightarrow}\limits_{\mbox{\tiny $x \!\rightarrow\! \infty$}}
0.
\end{equation}
Therefore, we get Theorem {\bf 3.6.} from \mbox{\bf \cite{Petojevic_06}}. Next we give a solution
to the open problem stated in Question {\bf 3.7.} in \mbox{\bf \cite{Petojevic_06}}. Namely,
the following formula in the article \mbox{\bf \cite[\mbox{\rm p.$\,$35.}]{Malesevic_03}} is given:
\begin{equation}
\label{Eq_10}
K(z)
=
\displaystyle\frac{\mbox{\rm Ei}(1)+\eu{i}\,\pi}{e}
+
\displaystyle\frac{(-1)^{z} \Gamma(1+z) \Gamma(-z,-1)}{e},
\end{equation}
for values $z \!\in\! \mbox{\bf C} \backslash \{-1,-2,-3,-4,\ldots\}$ and $\eu{i} \!=\! \sqrt{-1}$. In the previous
formula $\mbox{\rm Ei}(z)$ and $\Gamma(z,a)$ are exponential integral and incomplete gamma function respectively
\mbox{\bf \cite{Malesevic_03}}. Then, for fixed $i \!\in\! \mbox{\bf N}$ and values $z \!\in\! \mbox{\bf C}
\backslash \{-i,-i-1,-i-2,-i-3,\ldots\}$, on the basis of (\ref{Eq_4})~and~(\ref{Eq_10}), we~get:

\vspace*{-0.25 ex}

\noindent
\begin{equation}
\label{Eq_11}
\begin{array}{rcl}
K_{i}(z)
&\!\!=\!\!&
\displaystyle\frac{1}{(i\!-\!1)!}
{\Big (}K(z+i-1) - K(i-1){\Big )}                                              \\[2.0 ex]
&\!\!=\!\!&
\displaystyle\frac{\mbox{\rm Ei}(1) + \eu{i}\,\pi}{e (i\!-\!1)!}
+
\displaystyle\frac{(-1)^{z+i-1}\Gamma(1+z+i-1)\Gamma(-z-i+1,-1)}{e (i\!-\!1)!} \\[2.0 ex]
& &
\!\!\!\!
-
\displaystyle\frac{\mbox{\rm Ei}(1) + \eu{i}\,\pi}{e (i\!-\!1)!}
-
\displaystyle\frac{(-1)^{i-1}\Gamma(i)\Gamma(-i+1,-1)}{e (i\!-\!1)!}           \\[1.75 ex]
&\!\!=\!\!&
(-1)^{i}e^{-1}
{\bigg (}
\Gamma(1-i,-1) - (-1)^z \displaystyle\frac{\Gamma(1-i-z,-1)\Gamma(i+z)}{(i\!-\!1)!}
{\bigg )}.
\end{array}
\end{equation}
Therefore, the affirmative answer for Question {\bf 3.7.} from \mbox{\bf \cite{Petojevic_06}}
is true for complex values $z \!\in\! \mbox{\bf C} \backslash \{-i,-i-1,-i-2,-i-3, \ldots \}$.

\break

\medskip
\noindent
Finally, at the end of this note let us emphasize one differential--algebraic fact for the sequence
of functions $K_{i}(z)$. On the basis of the formula (17) from the article \mbox{\bf \cite{Petojevic_06}},
we can conclude that each $K_{i}(z)$ function satisfies the following recurrence relation
\mbox{$(i\!-\!1)!\,K_{i}(z+1)-(i\!-\!1)!\,K_{i}(z)=\Gamma(z+i)$}. The previous relation
can be used to verify the differential transcendency of these functions as discussed
in \mbox{\bf \cite{Mijajlovic_Malesevic_06, Mijajlovic_Malesevic_07}}. Therefore,
we can conclude that each $K_{i}(z)$ function is a differential transcendental function,
i.e. it satisfies no algebraic differential equation~over the~field~of~complex~rational~functions.

{\small
\noindent University of Belgrade,
          \hfill {\em $(\,$Received$\,:\;04/01/2007\;)$}              \break
\noindent Faculty of Electrical Engineering,
          \hfill {\em $(\,$Accepted$\,:\,05/25/2007\;)$}              \break
\noindent P.O.Box 35-54, $11 \, 120$ Belgrade, Serbia           \hfill\break
\noindent {\footnotesize {\bf malesh@eunet.yu}, {\bf malesevic@etf.bg.ac.yu}}
\hfill}

\break

\end{document}